\documentclass[11pt]{article}

\usepackage{amsmath}
\usepackage{amssymb}

 \newcommand{\beq}{\begin{equation}}
\newcommand{\eeq}{\end{equation}}

\numberwithin{equation}{section}

\newtheorem{theorem+}           {Theorem}      [section]
\newtheorem{definition+}  [theorem+]  {Definition}
\newtheorem{lemma+}  [theorem+]  {Lemma}
\newtheorem{corollary+}  [theorem+]  {Corollary}
\newtheorem{proposition+}  [theorem+]  {Proposition}
\newtheorem{example+}  [theorem+]  {Example}
\newtheorem{examples+}  [theorem+]  {Examples}
\newtheorem{remark+}  [theorem+]  {Remark}
\newtheorem{remarks+}  [theorem+]  {Remarks}

\newenvironment{theorem}{\begin{theorem+}\sl}{\end{theorem+}\rm}

\newenvironment{lemma}{\begin{lemma+}\sl}{\end{lemma+}\rm}
\newenvironment{corollary}{\begin{corollary+}\sl}{\end{corollary+}\rm}
\newenvironment{proposition}{\begin{proposition+}\sl}{\end{proposition+}\rm}

\newenvironment{examples}{\begin{examples+}\rm}{\end{examples+}\rm}
\newenvironment{proof}{\medbreak\noindent{\it Proof.}\rm}{\hfill$\square$\rm}
\newenvironment{remark}{\begin{remark+}\rm}{\end{remark+}\rm}
\newenvironment{remarks}{\begin{remarks+}\rm}{\end{remarks+}\rm}

\newcommand{\nux}{\nu_u(\zeta)}
\newcommand{\nuxa}{\nu_u(\zeta,a)}
\newcommand{\nuph}{\nu(u,\vph)}
\newcommand{\suph}{\sigma(u,\vph)}

\newcommand{\obl}{\Omega}

\newcommand{\okn}{1\le k\le n}
\newcommand{\vph}{\varphi}
\newcommand{\vphz}{\varphi_\zeta}

\renewcommand{\Bbb}{\mathbb}
\newcommand{\Z}{{\Bbb  Z}}
\newcommand{\Znp}{{ \Bbb Z}_+^n}

\newcommand{\Rnm}{{ \Bbb R}_-^n}
\newcommand{\Rnp}{{ \Bbb R}_+^n}
\newcommand{\C}{{\Bbb  C}}
\newcommand{\Cn}{{\Bbb  C\sp n}}

\newcommand{\D}{{\Bbb D}}

\newcommand{\F}{{\mathcal F}}

\newcommand{\M}{{\mathcal M}}

\newcommand{\I}{{\mathcal I}}
\newcommand{\Oo}{{\mathcal O}}

\newcommand{\supp}{{\operatorname{supp}\,}}
\newcommand{\codim}{{\operatorname{codim}\,}}

\begin{document}

\begin{center}
{\Large\bf Relative types and extremal problems\\ for
plurisubharmonic functions}
\end{center}

\medskip
\begin{center}
{\large\bf Alexander Rashkovskii}
\end{center}

\vskip3cm

\begin{abstract}\noindent
A type $\suph$ of a plurisubharmonic function $u$ relative to a
maximal plurisubharmonic weight $\vph$ with isolated singularity at
$\zeta$ is defined as $\liminf u(x)/\vph(x)$ as $x\to\zeta$. We
study properties of the relative types as functionals $u\mapsto
\suph$; it is shown that they give a general form for upper
semicontinuous, positive homogeneous and tropically additive
functionals on plurisubharmonic singularities. We consider some
extremal problems whose solutions are Green-like functions that give
best possible bounds on $u$, given the values of its types relative
to some of (or all) weights $\vph$; in certain cases they coincide
with known variants of pluricomplex Green functions. An analyticity
theorem is proved for the upperlevel sets for the types with respect
to exponentially H\"older continuous weights, which leads to a
result on propagation of plurisubharmonic singularities.

\end{abstract}

\section{Introduction}

If a holomorphic mapping $F$ vanishes at a point $\zeta$, then the
asymptotic behaviour of $|F|$ near $\zeta$ completely determines
such fundamental characteristics of $F$ at $\zeta$ as the
multiplicity of the zero or the integrability index. On the other
hand, in most cases the values of such characteristics can just give
certain bounds on the asymptotics of $F$ rather than recover it
completely.

The transformation $F\mapsto\log|F|$ puts this into the context of
pluripotential theory, which leads to a question of
characteristics of singularities of plurisubharmonic functions and
their relations to the asymptotic behaviour of the functions.
Since our considerations are local, we assume the functions to be
defined on domains of $\C^n$, $n>1$.

Let $u$ be a plurisubharmonic function near a point $\zeta$ of
$\Cn$, such that $u(\zeta)=-\infty$. The value $\nu_u(\zeta)$ of the
Lelong number of $u$ at $\zeta$ gives some information on the
asymptotic behaviour near $\zeta$: $u(x)\le
\nu_u(\zeta)\log|x-\zeta|+O(1)$. A more detailed information can be
obtained by means of its directional Lelong numbers
$\nu_u(\zeta,a)$, $a\in\Rnp$, due to Kiselman:
$u(x)\le\nu_u(\zeta,a)\max_k \,a_k^{-1}\log|x_k-\zeta_k|+O(1)$.

In addition, these characteristics of singularity are well suited
for the tropical structure of the cone of plurisubharmonic
functions, namely $\nu_{u+v}=\nu_u+\nu_v$ ({\it tropical
multiplicativity}) and
 $\nu_{\max\{u,v\}}=\min\{\nu_u,\nu_v\}$ ({\it tropical
additivity}). These properties play an important role, for example,
in investigation of valuations on germs of holomorphic functions
\cite{FaJ}. Note that the tropical operations $u\oplus
v:=\max\{u,v\}$ and $u\otimes v:=u+v$, when applied to
plurisubharmonic singularities, can be viewed as Maslov's
dequantization of usual addition and multiplication of holomorphic
functions.

A general notion of Lelong numbers $\nu(u,\vph)$ with respect to
plurisubharmonic weights $\vph$ was introduced and studied by
Demailly \cite{D1}, \cite{D}. Due to their flexibility, the
Lelong--Demailly numbers have become a powerful tool in
pluripotential theory and its applications.  They still are
tropically multiplicative, however tropical additivity is no longer
true for $\nu(u,\vph)$ with arbitrary plurisubharmonic weights
$\vph$, even if they are maximal outside $\vph^{-1}(-\infty)$. In
addition, the value $\nu(u,\vph)$ gives little information on the
asymptotics of $u$ near $\vph^{-1}(-\infty)$.

The good properties of the classical and directional Lelong numbers
result from the fact that they can be evaluated by means of the
suprema of $u$ over the corresponding domains (the balls and
polydiscs, respectively). This makes it reasonable to study the
asymptotics of the suprema of $u$ over the corresponding domains
$\{\vph(x)<t\}$ for a maximal weight $\vph$ with an isolated
singularity at $\zeta$ and consider the value $\suph=\liminf
u(x)/\vph(x)$ as $x\to\zeta$,  the {\it type of $u$ relative to}
$\vph$. The relative type is thus an alternative generalization of
the notion of Lelong number.

Unlike the Lelong--Demailly numbers, the relative types need not be
tropically multiplicative, however they are tropically additive.
Moreover, they are the only "reasonable" tropically additive
functionals on plurisubharmonic singularities (for a precise
statement, see Theorem~\ref{theo:repre}).

Maximality of $\vph$ gives the bound $u\le \suph\vph+O(1)$ near the
pole of $\vph$. We are then interested in best possible bounds on
$u$, given the values of its types relative to some of (or all) the
weights $\vph$ with fixed $\vph^{-1}(-\infty)$. Tropical additivity
of the relative types makes them a perfect tool for dealing with
upper envelopes of families of plurisubharmonic functions,
constructing thus extremal plurisubharmonic functions with
prescribed singularities. In certain cases these Green-like
functions coincide with known variants of pluricomplex Green
functions. In particular, this gives a new representation of the
Green functions with divisorial singularities
(Theorems~\ref{theo:grf1} and \ref{theo:grf2}). We study relations
between such extremal functions; one of the relations implies a
complete characterization of holomorphic mappings $f$ with isolated
zero at $\zeta$ of the multiplicity equal to the Newton number of
$f$ at $\zeta$ (Corollary~\ref{cor:newt}).

We also prove that the upperlevel sets for the types relative to
exponentially H\"older continuous weights are analytic varieties (an
analogue to the Siu theorem). As an application, we obtain a result
on propagation of plurisubharmonic singularities
(Corollary~\ref{cor:propag}) that results in a new representation of
the Green functions with singularities along complex spaces
(Corollary~\ref{cor:grf3}).

\medskip
The paper is organized as follows. Section~\ref{sec:prelim} recalls
basic facts on Lelong numbers and Green functions. In
Section~\ref{sec:types} we present the definition and elementary
properties of the relative types. A representation theorem for
tropically additive functionals on plurisubharmonic singularities is
proved in Section~\ref{sec:repr}. In Sections~\ref{sec:Green} and
\ref{sec:gr&gr} we consider extremal problems for plurisubharmonic
functions with given singularities. An analyticity theorem for the
upperlevel sets  and its applications are presented in
Section~\ref{sec:Hoelder}.

\medskip

\section{Preliminaries}\label{sec:prelim}

\subsection{Lelong numbers}

The {\it Lelong number} $\nu_T(\zeta)$ of a closed positive current
$T$ of bidimension $(p,p)$ at a point $\zeta\in\Cn$ is the residual
mass of $T\wedge (dd^c\log|\cdot-\zeta|)^{p}$ at $\zeta$:
\beq\label{eq:LNcl} \nu_T(\zeta)=\lim_{r\to 0}\int_{|x-\zeta|<r}
T\wedge (dd^c\log|x-\zeta|)^{p}; \eeq here $d=\partial +
\bar\partial,\ d^c= ( \partial -\bar\partial)/2\pi i$.

The Lelong number $\nux$ of a plurisubharmonic function $u$ is just
the Lelong number of the current $dd^cu$. It can also be calculated
as \beq \label{eq:form1}\nux=\lim_{r\to
-\infty}r^{-1}\int_{S_1}u(\zeta+xe^r)\,dS_1(x), \eeq where $dS_1$ is
the normalized Lebesgue measure on the unit sphere $S_1$, as well as
\beq \label{eq:form2}\nux=\lim_{r\to -\infty}r^{-1}\sup\{u(x):\:
|x-\zeta|< e^r\}=\liminf_{z\to \zeta}\frac{u(z)}{\log|z-\zeta|},
\eeq see \cite{Kis1}. Since the function $\sup\{u(x):\: |x-\zeta|<
e^r\}$ is convex in $r$, representation (\ref{eq:form2}) implies the
bound $u(x)\le \nux\log|x-\zeta|+O(1)$ near $\zeta$.

Lelong numbers are independent of the choice of coordinates. Siu's
theorem states that the set $\{\zeta:\:\nu_T(\zeta)\ge c\}$ is
analytic for any $c>0$.

\subsection{Directional Lelong numbers}
A more detailed information on the behaviour of $u$ near $\zeta$ can
be obtained by means of the {\it directional Lelong numbers} due to
Kiselman \cite{Kis2}: given $a=(a_1,\ldots,a_n)\in\Rnp$,
 \beq\label{eq:form4}\nuxa=\lim_{r\to
-\infty}r^{-1}\sup\{u(x):\: |x_k-\zeta_k|<e^{ra_k},\ \okn\},\eeq or
equivalently, in terms of the mean values of $u$ over the
distinguished boundaries of the polydiscs, similarly to
(\ref{eq:form1}). Namely,
$u(x)\le\nu_u(\zeta,a)\phi_{a,\zeta}(x)+O(1)$ with \beq
\phi_{a,\zeta}(x)=\max_k \,a_k^{-1}\log|x_k-\zeta_k|.
\label{eq:phax} \eeq

Analyticity of the upperlevel sets $\{\zeta:\:\nuxa\ge c\}$  was
established in \cite{Kis2}, \cite{Kis3}.

Directional Lelong numbers give rise to the notion of local
indicators of plurisubharmonic functions \cite{LeR}. Given a
plurisubharmonic function $u$, its {\it (local) indicator} at a
point $\zeta$ is a plurisubharmonic function $\Psi_{u,\zeta}$ in
the unit polydisc $\D^n$ such that for any $y\in\D^n$ with
$y_1\cdot\ldots\cdot y_n\neq 0$, \beq
\Psi_{u,\zeta}(y)=-\nu_u(\zeta,a),\quad
a=-(\log|y_1|,\ldots,\log|y_n|)\in\Rnp. \label{eq:lind} \eeq It is
the largest nonpositive plurisubharmonic function in $\D^n$ whose
directional Lelong numbers at $0$ coincide with those of $u$ at
$\zeta$,  see the details in \cite{LeR}, \cite{R}.

\medskip

\subsection{Lelong--Demailly numbers}

A general notion of Lelong numbers with respect to plurisubharmonic
weights was introduced and studied by J.-P.~Demailly \cite{D1},
\cite{D}. Let $T$ be a closed positive current of bidimension
$(p,p)$ on a domain $\obl\subset\Cn$, and let $\vph$ be a continuous
plurisubharmonic function $\obl\to [-\infty,\infty)$, semiexhaustive
on the support of $T$, that is, $B_R^\vph\cap \supp T\Subset\obl$
for some real $R$, where $B_R^\vph:=\{x:\: \vph(x)<R\}$, and let
$S_{-\infty}^\vph:=\vph^{-1}(-\infty)\neq \emptyset$. The value
\beq\nu(T,\vph)=\lim_{r\to -\infty}\int_{B_r^\vph} T\wedge(dd^c\vph
)^{p}= T\wedge(dd^c\vph )^{p}(S_{-\infty}^\vph)\eeq is called {\it
the generalized Lelong number}, or {\it the Lelong--Demailly
number}, of $T$ with respect to the weight $\vph$. When
$\vph(x)=\log|x-\zeta|$, this is just the classical Lelong number of
$T$ at $\zeta$.

For a plurisubharmonic function $u$, we use the notation
$\nuph=\nu(dd^cu,\vph)$.

The generalized Lelong numbers have the following semicontinuity
property.

\begin{theorem}\label{theo:Dsc} {\rm (\cite{D}, Prop.~3.11)} If
$T_k\to T$ and $\vph$ is semiexhaustive on a closed set containing
the supports of all $T_k$, then
$\limsup_{k\to\infty}\nu(T_k,\vph)\le \nu(T,\vph)$.
\end{theorem}

The following comparison theorems describe variation of the
Lelong--Demailly numbers with respect to the weights and currents.

\begin{theorem}\label{theo:CT1} {\rm (\cite{D}, Th. 5.1)} Let $T$ be
a closed positive current of bidimension $(p,p)$, and let $\vph$,
$\psi$ be two weights semiexhaustive on $\supp T$ such that
$\limsup\,\psi(x)/\vph(x)=l<\infty$ as $\vph(x)\to -\infty$. Then
$\nu(T,\psi)\le l^p\,\nu(T,\vph)$.
\end{theorem}

\begin{theorem}\label{theo:CT2} {\rm (\cite{D}, Th. 5.9)} Let $u$ and
$v$ be plurisubharmonic functions such that $(dd^cu)^q\wedge T$ and
$(dd^cv)^q\wedge T$ are well defined near $S_{-\infty}^\vph$ and
$u=-\infty$ on $\supp T\cap S_{-\infty}^\vph$, where $T$ is a closed
positive current of bidimension $(p,p)$, $p\ge q$. If
$\limsup\,v(x)/u(x)=l<\infty$ as $\vph(x)\to -\infty$, then
$\nu((dd^cv)^q\wedge T,\vph)\le l^q\,\nu((dd^cu)^q\wedge T,\vph)$.
\end{theorem}

The directional Lelong numbers can be expressed in terms of the
Lelong--Demailly numbers with respect to the directional weights
$\phi_{a,\zeta}$ (\ref{eq:phax}):
\beq\label{eq:dirKD}\nuxa=a_1\ldots a_n\,\nu(u,\phi_{a,\zeta}).\eeq

Siu's theorem on analyticity of upperlevel sets was extended  in
\cite{D1} to generalized Lelong numbers with respect to weights
$\vphz(x)=\vph(x,\zeta)$ that are exponentially H\"older continuous
with respect to $\zeta$.

\medskip

\subsection{Green functions}

Let $D$ be a hyperconvex domain in $\Cn$ and let $PSH^-(D)$ denote
the class of all negative plurisubharmonic functions in $D$.

\medskip
The {\it pluricomplex Green function $G_{\zeta,D}$ of $D$ with
logarithmic pole at} $\zeta\in D$ (introduced by Lempert, Zahariuta,
Klimek) is the upper envelope of the class $\F_{\zeta,D}$ of $u \in
PSH^-(D)$ such that $u(x)\leq \log|{x-\zeta}|+O(1)$ near $\zeta$.
The class $\F_{\zeta,D}$ can be also described as the collection of
$u\in PSH^-(D)$ such that $\nu_u(\zeta)\ge 1$. The function
satisfies $G_{\zeta,D}(x)=\log|x-\zeta|+O(1)$ near $\zeta$ and
$(dd^cG_{\zeta,D})^n=\delta_\zeta$.

\medskip
A more general construction was presented by Zahariuta \cite{Za0},
\cite{Za}. Given a continuous plurisubharmonic function $\vph$ in a
neighbourhood of $\zeta\in D$ such that $\vph^{-1}(-\infty)=\zeta$
and $(dd^c\vph)^n=0$ outside $\zeta$, let \beq\label{eq:gzf}
G_{\vph,D}(x)=\sup\{u(x):\: u\in PSH^-(D),\ u\leq \vph+O(1)\ {\rm
near\ }\zeta\}.\eeq Then $G_{\vph,D}$ is maximal in
$D\setminus\{\zeta\}$ and $G_{\vph,D}=\vph+O(1)$ near $\zeta$. We
will refer to this function as the {\it Green--Zahariuta function
with respect to the singularity} $\vph$.

\medskip
A Green function with prescribed values of all directional Lelong
numbers at $\zeta$ (the {\it Green function with respect to an
indicator} $\Psi$) was introduced in \cite{LeR} as
\beq\label{eq:LNind} G_{\Psi,\zeta,D}(x)=\sup\{u(x):\: u\in
PSH^-(D),\ \nu_u(\zeta,a)\geq \nu_\Psi(0,a)\quad \forall
a\in\Rnp\},\eeq where $\Psi$ is a negative plurisubharmonic function
in the unit polydisc $\D^n$ such that \beq \label{eq:defind}
\Psi(z_1,\ldots,z_n)=\Psi(|z_1|,\ldots,|z_n|)=c^{-1}\Psi(|z_1|^c,
\ldots,|z_n|^c),\quad c>0,\ z\in\D^n;\eeq such a function $\Psi$
coincides with its own indicator (\ref{eq:lind}), so
$\nu_\Psi(0,a)=-\Psi(e^{-a_1},\ldots,e^{-a_n})$.

\medskip
The above Green functions were also considered for several isolated
poles. In the case of non-isolated singularities, a variant of Green
functions was introduced by L\'arusson--Sigurdsson \cite{LarSig1} by
means of the class of negative plurisubharmonic functions $u$
satisfying $\nu_u(x)\ge \alpha(x)$, where $\alpha$ is an arbitrary
nonnegative function on $D$. In \cite{LarSig2} it was specified to
the case when $\alpha(x)=\nu_A(x)$ is the Lelong number of a divisor
$A$; then \beq\label{eq:LSG2} \F_{A,D}=\{u\in PSH^-(D):\:
\nu_u(x)\geq \nu_A(x)\quad \forall x\in D\},\eeq and
\beq\label{eq:LSG1} G_{A,D}(x)=\sup\{u(x):\: u\in \F_{A,D}\}\eeq  is
the {\it Green function for the divisor} $A$. It was shown that if
$A$ is the divisor of a bounded holomorphic function $f$ in $D$,
then $G_{A,D}=\log|f|+O(1)$ near points of $|A|=f^{-1}(0)$.

This was used in \cite{RSig1}, \cite{RSig2} for consideration of
Green functions with arbitrary analytic singularities. Given a
closed complex subspace $A$ on $D$, let
${\I}_A=\I_{A,D}=({\I}_{A,x})_{x\in D}$ be the associated coherent
sheaf of ideals in the sheaf $\Oo_D=(\Oo_x)_{x\in D}$ of germs of
holomorphic functions on $D$, and let $|A|=\{x\in D:\:
{\I}_{A,x}\neq \Oo_x\}$. A {\it Green function $G_{A,D}$ for the
complex space} $A$ on $D$ was constructed as \beq\label{eq:RS}
G_{A,D}(x)=\sup\,\{u(x):\: u\in PSH^-(D), \ u\le \log|f|+O(1)\ {\rm
locally\ near \ }|A|\},\eeq where $f=(f_1,\ldots,f_p)$ and
$f_1,\ldots,f_p$ are local generators of $\I_A$. The function
$G_{A,D}$ is plurisubharmonic and satisfies \beq\label{eq:RS1}
G_{A,D}\le\log|f|+O(1)\eeq locally near points of $|A|$; if $\I_A$
has bounded global generators, then $G_{A,D}=\log|f|+O(1)$.


\medskip

\section{Definition and elementary properties of relative
types}\label{sec:types}

Given $\zeta\in\Cn$, let $PSH_\zeta$ stand for the collection of all
(germs of) plurisubharmonic functions $u\not\equiv-\infty$ in a
neighbourhood of $\zeta$.

Let $\vph\in PSH_\zeta$ be locally bounded outside $\zeta $,
$\vph(\zeta)=-\infty$, and maximal in a punctured neighbourhood of
$\zeta$: $(dd^c\vph)^n=0$ on $B_R^\vph\setminus\{\zeta\}$ for some
$R>-\infty$; we recall that $B_R^\vph=\{z:\: \vph(z)<R\}$. The
collection of all such {\it maximal weights} (centered at $\zeta$)
will be denoted by $MW_\zeta$. If we want to specify that
$(dd^c\vph)^n=0$ on $\omega\setminus\{\zeta\}$, we will write
$\vph\in MW_\zeta(\omega)$.

Given a function $u\in PSH_\zeta$, its singularity at $\zeta$ can be
compared to that of $\vph\in MW_\zeta$; the value
\beq\label{eq:type} \sigma(u,\vph)=\liminf_{z\to
\zeta}\frac{u(z)}{\vph(z)}\eeq  will be called the $\vph$-{\it
type}, or the {\it relative type} of $u$ with respect to $\vph$.

Both the Lelong--Demailly numbers and the relative types are
generalizations of the classical notion of Lelong number, however
they use different points of view on the Lelong number: while the
Lelong--Demailly numbers correspond to (\ref{eq:LNcl}) (and
(\ref{eq:form1}) in the case of functions), the relative types are
based on (\ref{eq:form2}). As we will see, the two generalizations
have much in common, however some features are quite different.

Note that $\sigma(u+h,\vph)=\sigma(u,\vph)$ if $h$ is
pluriharmonic (and hence bounded) near $\zeta$, so $\suph$ is
actually a function of $dd^cu$.

The following properties are direct consequences of the definition
of the  relative type.

\begin{proposition} \label{prop:rel} Let $\vph, \psi\in MW_\zeta$ and $u,
u_j\in PSH_\zeta$. Then
\begin{enumerate} \item[(i)]
$\sigma(cu,\vph)=c\,\suph$ for all $c>0$; \item[(ii)]
$\sigma(\max\{u_1,u_2\},\vph)=\min_j\sigma(u_j,\vph)$;
\item[(iii)]
$\sigma(u_1+u_2,\vph)\ge\sigma(u_1,\vph)+\sigma(u_2,\vph)$;
\item[(iv)]$\sigma(u,\psi)\ge \suph\,\sigma(\vph,\psi)$; in
particular, if there exists
$\lim_{z\to\zeta}\vph(z)/\psi(z)=l<\infty$, then $\suph=
l\,\sigma(u,\psi)$; \item[(v)] if
$\liminf_{z\to\zeta}u_1(z)/u_2(z)=l<\infty$, then
$\sigma(u_1,\vph)\ge l\,\sigma(u_2,\vph)$; \item[(vi)] if
$u=\log\sum_j^m e^{u_j}$, then
$\sigma(u,\vph)=\sigma(\max_ju_j,\vph)=\min_j\sigma(u_j,\vph)$.
\end{enumerate}
\end{proposition}

\begin{proof} Properties (i) -- (v) are direct consequences of the definition
of the  relative type, and (vi) follows from (ii) and (iv) together.
Note that (iv) makes sense because $\suph<\infty$ for any $u\in
PSH_\zeta$ and $\vph\in MW_\zeta$; this follows, for example, from
an alternative description of $\suph$ given by
(\ref{eq:luvr})--(\ref{eq:ratlim}) below.
\end{proof}


\begin{remarks}\label{rem:3.2} (1) The relative types need not be tropically
multiplicative functionals on $PSH_\zeta$. Take
$\vph=\max\,\{3\log|z_1|,3\log|z_2|, \log|z_1z_2|\}\in MW_0$ and
$u_j=\log|z_j|$, $j=1,2$,  then $\sigma(u_j,\vph)=1/3$, while
$\sigma(u_1+u_2,\vph)=1$.

(2) Properties (iv) and (v) are analogues to Comparison Theorems
\ref{theo:CT1} and \ref{theo:CT2}.

(3) For holomorphic functions $f_j$ and positive numbers $p_j$,
$j=1,\ldots,m$, property (vi) gives the relation
\beq\sigma(\log\sum_j |f_j|^{p_j},\vph)= \sigma(\max_j
p_j\log|f_j|,\vph)=\min_j p_j\,\sigma(\log|f_j|,\vph).\eeq
\end{remarks}

Given a weight $\vph\in MW_\zeta$ and a function $u\in PSH_\zeta$,
consider the growth function \beq\label{eq:luvr}
\Lambda(u,\vph,r)=\sup\,\{u(z):\:z\in B_r^\vph\}.\eeq

\begin{proposition} \label{prop:convexity}
{\rm (see also \cite{D3}, Corollary 6.6)} Let $\vph\in
MW_\zeta(B_R^\vph)$ such that $B_R^\vph$ is bounded, and $u\in
PSH(B_R^\vph)$. Then the function $\Lambda(u,\vph,r)$ is convex in
$r\in (-\infty,R)$.
\end{proposition}

\begin{proof} Given $-\infty<r_1<r_2< R$ and $0<\epsilon <R-r_2$, let
$u_\epsilon=u*\chi_\epsilon$ be a standard regularization
(smoothing) of $u$. Take $c\ge 0$ and $d\in{\Bbb R}$ such that
$cr_j+d=\Lambda(u_\epsilon,\vph,r_j)$, $j=1,2$. Since $\vph\ge r$ on
$\partial B_{r}^\vph$, we have $c\vph+d\ge u_\epsilon$ on
$\partial(B_{r_2}^\vph\setminus \overline{B_{r_1}^\vph})$ and thus
on the set $B_{r_2}^\vph\setminus B_{r_1}^\vph$ because of the
maximality of the function $c\vph+d$ on $B_R^\vph\setminus\{0\}$. In
other words, $\Lambda(u_\epsilon,\vph,r)\le cr+d$ on $(r_1,r_2)$,
which means that $\Lambda(u_\epsilon,\vph,r)$ is convex on
$(r_1,r_2)$. Since $\Lambda(u_\epsilon,\vph,r)\to \Lambda(u,\vph,r)$
as $\epsilon\to 0$, this implies the assertion.
\end{proof}

\medskip

Since $\Lambda(u,\vph,r)$ is increasing and convex, the ratio \beq
g(u,\vph,r,r_0):=\frac{\Lambda(u,\vph,r)-\Lambda(u,\vph,r_0)}{r-r_0},\quad
r<r_0< R,\eeq is increasing in $r\in(-\infty,r_0)$ and, therefore,
has a limit as $r\to-\infty$; it is easy to see that the limit
equals $\suph$. We have, in particular, \beq\label{eq:ratlim}
\sigma(u,\vph)\le g(u,\vph,r,r_0),\quad r<r_0,\eeq which implies the
following basic bound.

\begin{proposition} \label{prop:phbound} Let $\vph\in
MW_\zeta(B_R^\vph)$, $u\in PSH^-(B_{r_0}^\vph)$, $r_0< R$. Then
$u\le\suph (\vph-r_0)$ in $B_{r_0}^\vph$. In particular, every
function $u\in PSH_\zeta$ has the bound
 \beq\label{eq:lb}
u(z)\le\suph\,\vph(z)+O(1),\quad z\to \zeta.\eeq
\end{proposition}

\medskip
Next statement is an analogue to Theorem~\ref{theo:Dsc}.

\begin{proposition} \label{prop:converg} Let $u_j,u\in
PSH^-(\obl)$, $\vph\in MW_\zeta$, $\zeta\in\obl$. If $u_j\to u$ in
$L_{loc}^1(\obl)$, then $\suph\ge\limsup\,\sigma(u_j,\vph).$
\end{proposition}

\begin{proof} Take any $r_0$ such that $\Lambda(u,\vph,r_0)<0$,
then for any $\epsilon>0$ and $r<r_0$ there exists $j_0$ such that
$\Lambda(u_j,\vph,r)\ge \Lambda(u,\vph,r)-\epsilon $ for all $j>
j_0$. Using (\ref{eq:ratlim}) we get \beq \sigma(u_j,\vph)\le
g(u_j,\vph,r,r_0)\le \frac{\Lambda(u,\vph,r)-\epsilon}{r-r_0},\eeq
which implies the assertion. \end{proof}

\medskip

Let us compare the values of relative types with some known
characteristics of plurisubharmonic singularities. Denote
\beq\nu_\vph:=\nu_\vph(\zeta)\eeq the Lelong number of $\vph\in
MW_\zeta$ at $\zeta$; \beq\tau_\vph:=(dd^c\vph)^n(\zeta)\eeq the
residual Monge--Amp\`ere mass of $\vph$ at $\zeta$;
\beq\label{eq:uln}\alpha_\vph:=\limsup_{z\to
\zeta}\frac{\vph(z)}{\log|z-\zeta|}.\eeq By Theorem~\ref{theo:CT1},
$\nu_\vph^n\le\tau_\vph\le \alpha_\vph^n$ and, by
Proposition~\ref{prop:phbound}, \beq\label{eq:comparable}
\alpha_\vph\log|z-\zeta|+O(1)\le\vph(z)\le
\nu_\vph\log|z-\zeta|+O(1), \quad z\to \zeta.\eeq

If $\vph$ has analytic singularity, that is, $\vph=\log|f|+O(1)$
near $\zeta$, where $f=(f_1,\ldots,f_n)$ is a holomorphic map with
isolated zero at $\zeta$, then $\nu_\vph$ equals the minimum of the
multiplicities of $f_k$ at $\zeta$, $\tau_\vph$ is the multiplicity
of $f$ at $\zeta$, and $\alpha_\vph=\gamma_f$, the \L ojasiewicz
exponent of $f$ at $\zeta$, i.e., the infimum of $\gamma>0$ such
that $|f(z)|\ge |z-\zeta|^\gamma$ near $\zeta$. Therefore,
$\nu_\vph>0$ and $\alpha_\vph<\infty$ in this case.

In the general situation, since $\vph$ is locally bounded and
maximal on $B_r^\vph\setminus\{\zeta\}$, the condition
$\vph(\zeta)=-\infty$ implies $0<\tau_\vph<\infty$.

We do not know if $\nu_\vph>0$ for every weight $\vph\in MW_\zeta$.
It is actually equivalent to the famous problem of existence of a
plurisubharmonic function which is locally bounded outside $\zeta$
and has zero Lelong number and positive Monge-Amp\`ere mass at
$\zeta$ (see a discussion in \cite{W} or the remark after
Proposition~\ref{prop:geqh}).

Furthermore, (\ref{eq:comparable}) implies $\alpha_\vph>0$, however
we do not know if the "\L ojasiewicz exponent" $\alpha_\vph$ is
finite for every maximal weight $\vph$. It is worth noting that, by
Proposition~\ref{prop:rel} (v), $\sigma(\log|f|,\vph)\le
\gamma_f\alpha_\vph^{-1}$ for any holomorphic map $f$ with isolated
zero at $\zeta$, where $\gamma_f$ is the \L ojasiewicz exponent of
$f$, so $\sigma(\log|f|,\vph)=0$ if $\alpha_\vph=\infty$.

By Theorem~\ref{theo:CT2}, the condition $\alpha_\vph<\infty$
implies $\tau_\vph\le\alpha_\vph^{n-1}\nu_\vph$ and so, $\nu_\vph>0$
for such a weight $\vph$. In other words, denote \beq\label{eq:SMW}
SMW_\zeta=\{\vph\in MW_\zeta:\: \nu_\vph>0\}\eeq (the weights with
"strong" singularity) and \beq\label{eq:FMW} LMW_\zeta=\{\vph\in
MW_\zeta:\: \alpha_\vph<\infty\}\eeq (the weights with finite \L
ojasiewicz exponent), then \beq\label{eq:sb} LMW_\zeta\subset
SMW_\zeta\subset MW_\zeta,\eeq and it is unclear if the inclusions
are strict.


\begin{proposition} \label{prop:relsnu}
The type $\suph$ of $u\in PSH_\zeta$ with respect to $\vph\in
SMW_\zeta$ is related to the Lelong number $\nu_u(\zeta)$ of $u$ at
$\zeta$ as \beq\label{eq:lelbound}
\alpha_\vph^{-1}\nu_u(\zeta)\le\suph\le\nu_\vph^{-1}\nu_u(\zeta).\eeq
If, in addition, $\vph$ is continuous, then
\beq\label{eq:sunu}\suph\le\tau_\vph^{-1}\nuph,\eeq where $\nuph$ is
the Lelong--Demailly number of $u$ with respect to $\vph$.
\end{proposition}

\begin{proof} Bounds (\ref{eq:lelbound}) follow from
(\ref{eq:comparable}) by Theorem~\ref{theo:CT1}, and relation
(\ref{eq:sunu}) follows from Proposition~\ref{prop:phbound} by
Theorem~\ref{theo:CT2}.
\end{proof}


\begin{remark} Due to (\ref{eq:form4}) and (\ref{eq:dirKD}), there is
always an equality in (\ref{eq:sunu}) if $\vph=\phi_{a,\zeta}$, a
directional weight (\ref{eq:phax}). On the other hand, for general
weights $\vph$ the inequality can be strict. For example, let $u_1$,
$u_2$ and $\vph$ be as in Remarks~\ref{rem:3.2} (1), and let
$u=\max\,\{2u_1,u_2\}$. Then $\suph=1/3$, $\nuph=3$ and
$\tau_\vph=6$, so the right hand side of (\ref{eq:sunu}) equals
$1/2>1/3$.
\end{remark}



\section{Representation theorem}\label{sec:repr}

It was shown in Section~\ref{sec:types} that relative types $\suph$
are positive homogeneous, tropically additive (in the sense
$\sigma(\max u_k,\vph)=\min\, \sigma(u_k,\vph)$) and upper
semicontinuous functionals on $PSH_\zeta$ that preserve ordering of
the singularities, i.e., $u\le v+O(1)$ implies $\suph\ge
\sigma(v,\vph)$.

Here we show that any such functional on $PSH_\zeta$ can be
represented as a relative type, provided it does not vanish on a
function that is locally bounded outside $\zeta$.

\begin{lemma}\label{lemma:repre} Let $D$ be a bounded hyperconvex
neighbourhood of a point $\zeta$, and let a function $\sigma:\:
PSH^-(D)\to [0,\infty]$ be such that $\sigma(u)< \infty$ if
$u\not\equiv-\infty$ and
\begin{enumerate} \item[(i)] $\sigma(cu)=c\,\sigma(u)$ for all $c>0$;
\item[(ii)] if $u_1\le u_2+O(1)$ near $\zeta$, then $\sigma(u_1)\ge \sigma(u_2)$;
\item[(iii)] $\sigma(\max_k u_k)=\min_k \sigma(u_k)$, $k=1,2$;
\item[(iv)] if $u_j\to u$ in $L_{loc}^1$, then $\limsup\,
\sigma(u_j)\le\sigma(u)$;
\item[(v)] $\sigma(w_0)>0$ for at least one  $w_0\in PSH^-(D)\cap L_{loc}^\infty
(D\setminus\zeta)$.
\end{enumerate}
Then there exists a unique function $\vph\in MW_\zeta(D)$,
$\vph(z)\to 0$ as $z\to\partial D$, such that
\beq\label{eq:rpr}\sigma(u)=\suph\quad\forall u\in PSH^-(D).\eeq If,
in addition, (v) is true with $w_0(z)=\log|z-\zeta|+O(1)$, then
$\vph\in LMW_\zeta$.
\end{lemma}

\begin{proof}
Denote $\vph(z)=\sup\,\{u(z):\: u\in\M\}$, where $\M=\{u\in
PSH^-(D):\: \sigma(u)\ge 1\}$.

By the Choquet lemma, there exists a sequence $u_j\in\M$ increasing
to a function $v$ such that $v^*=\vph^*\in PSH^-(D)$. Properties
(iii) and (iv) imply $v^*\in\M$, so $v^*\le\vph$. Therefore,
$\vph=\vph^*\in \M$. Evidently, $\sigma(\vph)=1$.

If $v\in PSH^-(D)$ satisfies $v\le\vph$ outside $\omega\Subset
D\setminus\{\zeta\}$, then $v\in \cal M$ and so, $v\le\vph$ in $D$.
Therefore, the function $\vph$ is maximal on $D\setminus\{\zeta\}$.
Furthermore, $\vph\in L_{loc}^\infty(D\setminus\{\zeta\})$ because
$\vph\ge w_0/\sigma(w_0)$. It is not hard to see that
$\vph(\zeta)=-\infty$. Indeed, assuming $\vph(\zeta)=A>-\infty$, the
maximality of $\vph$ on $\{\vph(z)<A\}$ gives $\vph\ge A$
everywhere, which contradicts $\sigma(\vph)>0$ in view of (ii). So,
$\vph\in MW_\zeta(D)$.

Standard arguments involving a negative exhaustion function of $D$
show that $\vph(z)\to 0$ as $z\to\partial D$.

By Proposition~\ref{prop:phbound}, $u\le\suph \vph+O(1)$ and thus,
by (i) and (ii), $\sigma(u)\ge \suph $ for every $u\in PSH^-(D)$.
This gives, in particular, $\suph=0$ if $\sigma(u)=0$. Let
$\sigma(u)>0$, then $u/\sigma(u)\in\cal M$, so $u\le \sigma(u)\vph$
and consequently, $\sigma(u,\vph)\ge \sigma(u)$. This proves
(\ref{eq:rpr}).

If $\psi$ is another weight from $MW_\zeta(D)$ with zero boundary
values on $\partial D$, representing the functional $\sigma$, then
$\psi\le\vph$. On the other hand, the relation $\sigma(\vph,\psi)=1$
implies that for any $\epsilon\in (0,1)$ we have $\vph\le
(1-\epsilon)\psi+\epsilon$ on a neighbourhood of $\zeta$ and near
$\partial D$ and thus, by the maximality of $\psi$ on
$D\setminus\{\zeta\}$, everywhere in $D$.

Finally, the last assertion follows from the relation $\vph\ge
w_0/\sigma(w_0)\in LMW_\zeta$.
\end{proof}

\begin{remarks} (1) Note that for the functional
$\sigma(u)=\sigma(u,\phi)$ with a {\sl continuous} weight $\phi\in
MW_\zeta$, the function $\vph$ constructed in the proof of
Lemma~\ref{lemma:repre} is just the Green--Zahariuta function for
the (continuous) singularity $\phi$ (\ref{eq:gzf}). We will keep
this name for the case of Green functions with respect to arbitrary
singularities $\phi\in MW_\zeta$, \beq\label{eq:gzf1}
G_{\phi,D}(z)=\sup\{u(z)\in PSH^-(D):\: u\leq \phi+O(1)\ {\rm near\
}\zeta\}.\eeq We have thus $G_{\phi,D}\in MW_\zeta(D)$,
$G_{\phi,D}=\phi+O(1)$ because $\sigma(G_{\phi,D},\phi)=1$, and
$G_{\phi,D}(z)\to 0$ as $z\to\partial D$ if $D$ is hyperconvex.

(2) Let $\preceq$ be a natural partial ordering on $MW_\zeta$: $\phi
\preceq\psi$ if $\sigma(u,\phi)\le\sigma(u,\psi)$ for any $u\in
PSH_\zeta$. It is easy to see that $\phi \preceq\psi\Leftrightarrow
\sigma(\phi,\psi)\ge 1 \Leftrightarrow G_{\phi,D}\le G_{\psi,D}$ for
some (and, consequently, for any) hyperconvex neighbourhood $D$ of
$\zeta$.
\end{remarks}

The following representation theorem is an easy consequence of
Lemma~\ref{lemma:repre}.

\begin{theorem}\label{theo:repre} Let a function $\sigma:\: PSH_\zeta\to
[0,\infty)$ satisfy conditions [i)--(v) of Lemma~\ref{lemma:repre}
for $u,u_k\in PSH_\zeta$ and $D$ a bounded hyperconvex neighbourhood
of $\zeta$. Then there exists a weight $\vph\in MW_\zeta$ such that
$\sigma(u)=\suph$ for every $u\in PSH_\zeta$. The representation is
essentially unique: if two weights $\vph$ and $\psi$ represent
$\sigma$, then $\vph=\psi+O(1)$ near $\zeta$. If, in addition, (v)
is true with $w_0(z)=\log|z-\zeta|$, then $\vph\in LMW_\zeta$.
\end{theorem}

\begin{proof} We may assume $\sigma(w_0)=1$. Let $\vph\in MW_\zeta$ be
the function constructed in Lemma~\ref{lemma:repre}. Exactly as in
the proof of the lemma, we get $\sigma(u)\ge \suph $ for every $u\in
PSH_\zeta$. To prove the reverse inequality, take any $u\in
PSH_\zeta$. The function $v_0=\max\,\{u,\sigma(u)w_0\}$ can be
extended from a neighbourhood of $\zeta$ to a plurisubharmonic
function $v$ on a neighbourhood of $\overline D$; by (iii),
$\sigma(u)=\sigma(v_0)=\sigma(v-\sup_Dv)=\sigma(v-\sup_Dv,\vph)=
\sigma(v_0,\vph)\le \suph$, so $\sigma(u)=\suph$.

If $\psi\in MW_\zeta$ is another weight representing the functional
$\sigma$, then $\sigma(\vph,\psi)=\sigma(\psi,\vph)=1$ and
$\psi=\vph+O(1)$ by Proposition~\ref{prop:phbound}.
\end{proof}

\begin{remark} Recall that a {\it valuation} on the local ring
$R_\zeta$ of germs of analytic functions $f$ at $\zeta$ is a
nonconstant function $\mu:\: R_\zeta\to [0,+\infty]$ such that
\beq\mu(f_1f_2)=\mu(f_1)+\mu(f_2),
\quad\mu(f_1+f_2)\ge\min\,\{\mu(f_1),\mu(f_2)\}, \quad \mu(1)=0;\eeq
a valuation $\mu$ is {\it centered} if $\mu(f)>0$ for every $f$ from
the maximal ideal $\mathfrak{m}_\zeta$, and {\it normalized} if
$\min\, \{\mu(f):f\in\mathfrak{m}_\zeta\}=1$. Every weight $\vph\in
MW_\zeta$ generates a functional $\sigma_\vph$ on $R_\zeta$,
$\sigma_\vph(f)=\sigma(\log|f|,\vph)$, with the properties
\beq\sigma_\vph(f_1f_2)\ge \sigma_\vph(f_1)+\sigma_\vph(f_2), \quad
\sigma_\vph(f_1+f_2)\ge\min\,\{\sigma_\vph(f_1),\sigma_\vph(f_2)\},
\quad\sigma_\vph(1)=0.\eeq Such a functional is thus a valuation if
the weight $\vph$ satisfies the additional condition
\beq\label{eq:trm}\sigma(u+v,\vph)=\sigma(u,\vph)+\sigma(v,\vph)\eeq
for any $u,v\in PSH_\zeta$ -- in other words, if
$u\mapsto\sigma(u,\vph)$ is tropically linear (both additive and
multiplicative);  $\sigma_\vph$ is centered if and only if $\vph\in
LMW_\zeta$, and normalized iff $\alpha_\vph=1$.

The weights $\phi_{a,\zeta}$ (\ref{eq:phax}) satisfy (\ref{eq:trm}),
and the corresponding functionals $\sigma_{\phi_{a,\zeta}}$ are
monomial valuations on $R_\zeta$; they are normalized, provided
$\min_ka_k=1$. It was shown in \cite{FaJ} that an important class of
valuations in $\C^2$ (quasimonomial valuations) can be realized as
$\sigma_\vph$ with certain weights $\vph\in LMW_\zeta$ satisfying
(\ref{eq:trm}) and $\alpha_\vph=1$, and all other normalized
valuations in $\C^2$ can be realized as limits of increasing
sequences of the quasimonomial ones. We believe the relative types
with respect to weights satisfying (\ref{eq:trm}) can be used in
investigation of valuations in higher dimensions.
\end{remark}



\section{Greenifications}\label{sec:Green}

In this section we consider some extremal problems for
plurisubharmonic functions with singularities determined by a given
plurisubharmonic function $u$. Solutions to these problems resemble
various Green functions mentioned in Section~\ref{sec:prelim} (and
in some cases just coincide with them), and we will call them {\it
greenifications} of the function $u$. This reflects the point of
view on Green functions as largest negative plurisubharmonic
functions with given singularities; different types of the Green
functions arise from different ways of measuring the singularities
(or different portions of information on the singularities used).
Note that the tropical additivity makes relative types an adequate
tool in constructing extremal plurisubharmonic functions as upper
envelopes.


\subsection{Type-greenifications}

Let a bounded domain $D$ contain a point $\zeta$ and let $u\in
PSH_\zeta$. Given a collection $P$ of weights $\phi\in MW_\zeta$,
denote \beq\label{eq:mu1} \M_u^P=\M^P_{u,\zeta,D}=\{v(x):\: v\in
PSH^-(D),\ \sigma(v,\phi) \ge \sigma(u,\phi)\quad \forall\phi\in
P\}\eeq and define the function \beq\label{eq:genind0}
h_{u}^P(x)=h^P_{u,\zeta,D}(x)=\sup\,\{v(x):\:v\in
\M^P_{u,\zeta,D}\}.\eeq We will write simply $\M_{u,\zeta,D}$ and
$h_{u,\zeta,D}$ if $P=MW_\zeta$. The function $h^P_{u,\zeta,D}$ will
be called the {\it type-greenification of $u$ with respect to the
collection} $P$, and $h_{u,\zeta,D}$ will be called just the {\it
type-greenification of $u$ at} $\zeta$.

Consideration of the functions $h_u^P$ with $P\neq MW_\zeta$ can be
useful in situations where the only information on the singularity
of $u$ available is the values of $\suph$ for certain selected
weights $\vph$. One more reason is that for some collections $P$,
the functions $h_u^P$ are quite easy to compute and, at the same
time, they can give a reasonably good information on the asymptotic
behaviour of $u$ (see Examples~2 and 3 after
Proposition~\ref{prop:ubind0}). Note that $h^P_{u}\ge h^Q_{u}\ge
h_{u}$ if $P\subset Q\subset MW_\zeta$.

\begin{proposition} \label{prop:genind0} Let $u\in PSH(D)$ be bounded
above in $D\ni\zeta$, and $P\subseteq MW_\zeta$. Then
\begin{enumerate} \item[(i)] $h^P_{u,\zeta,D}\in PSH^-(D)$;
\item[(ii)] $u\le h^P_{u,\zeta,D}+\sup_Du$; \item[(iii)]
$\sigma(u,\phi)= \sigma(h^P_{u,\zeta,D},\phi)$ for any weight
$\phi\in P$; \item[(iv)] $h^P_{u,\zeta,D}$ is maximal on
$D\setminus\{\zeta\}$; \item[(v)] if $D$ has a strong
plurisubharmonic barrier at a point $z\in
\partial D$ (i.e., if there exists $v\in PSH(D)$ such that
$\lim_{x\to z}v(x)=0$ and $\sup_{D\setminus U}v<0$ for every
neighbourhood $U$ of $z$) and if $u$ is bounded below near $z$,
then $h^P_{u,\zeta,D}(x)\to 0$ as ${x\to z}$.
\end{enumerate}
\end{proposition}

\begin{proof} Let $w$ denote the right hand side of
(\ref{eq:genind0}), then its upper regularization $w^*$ is
plurisubharmonic in $D$. By the Choquet lemma, there exists a
sequence $v_j\in \M_{u}^P$ such that $w^*=(\sup_jv_j)^*$. Denote
$w_k=\sup_{j\le k}v_j$. By Proposition~\ref{prop:rel}(ii),
$w_k\in\M_{u}^P$. Since the functions $w_k$ converge weakly to
$w^*$, Proposition~\ref{prop:converg} implies then $\sigma(w^*,\phi)
\ge \sigma(w_k,\phi)\ge\sigma(u,\phi)$ for any weight $\phi\in P$
and so, $w^*\in\M_{u}^P$, which proves (i) and gives, at the same
time, the inequality $\sigma(u,\phi)\le \sigma(h_{u}^P,\phi)$. The
reverse inequality, even for arbitrary weights $\phi\in MW_\zeta$,
follows from the (evident) assertion (ii) and completes the proof of
(iii).

If a function $v\in PSH(D)$ satisfies $v\le h^P_{u,\zeta,D}$ on
$D\setminus\omega$ for some open set $\omega \Subset
D\setminus\zeta$, then $\max\,\{v,h^P_{u,\zeta,D}\}\in\M_u^P$ and
therefore $v\le h^P_{u,\zeta,D}$ on $\omega$, which proves (iv).

Finally, to prove (v), take a neighbourhood $U$ of $z\in\partial D$
such that $u>t>-\infty$ on $U\cap D$ and choose $c>0$ such that
$v<t/c$ on $D\setminus U$, so $u>cv$ on $\partial U\cap D$. Let \beq
w(x)=\left\{\begin{array}{ll}\max\{u(x), cv(x)\}, &
\mbox{$x\in D\cap U$}\\
u(x), & \mbox{$x\in D\setminus U$,}\end{array}\right. \eeq then
$w\in\M^P_{u,\zeta,D}(u)$ and $\lim_{x\to z}w(x)=0$.
\end{proof}

\medskip
If $v\in PSH^-(D)$, then the relation $\sigma(v,\vph) \ge
\sigma(u,\vph)$ is equivalent to $v\le \sigma(u,\vph)G_{\vph,D}$.
Therefore, these Green-like functions $h^P_{u,\zeta,D}$ can be
described by means of the Green--Zahariuta functions $G_{\vph,D}$
(\ref{eq:gzf1}) as follows.

\begin{proposition} \label{prop:ubind0}
The function $h^P_{u,\zeta,D}$, $P\subset MW_\zeta$, is the largest
plurisubharmonic minorant of the family $\{\sigma(u,\vph)
G_{\vph,D}:\: \vph\in P\}$. In particular, if $\vph\in MW_\zeta$,
then $h_{\vph,\zeta,D}=h^\vph_{\vph,\zeta,D}=G_{\vph,D}$.
\end{proposition}

\begin{examples} (1) If $P$ consists of a single weight $\vph$, then
$h^P_{u,\zeta,D}=\sigma(u,\vph)\,G_{\vph,D}$. In particular, if
$\vph(x)=\log|x-\zeta|$, then
$h^P_{u,\zeta,D}=\nu_u(\zeta)\,G_{\zeta,D}$.

(2) Let $A$ be a finite subset of $\Rnp$ and let $P$ be the
collection of the weights $\phi_a=\phi_{a,0}$ (\ref{eq:phax}) with
$a\in A$. According to Proposition~\ref{prop:ubind0} and
(\ref{eq:form4}), the function $h^P_{u,0,D}$ is the largest
plurisubharmonic minorant of the family $\{\nu_u(0,a)
G_{\phi_{a},D}:\: a\in A\} $. Using methods from \cite{LeR} and
\cite{R}, it can then be shown that the minorant is the
Green--Zahariuta function  for the singularity
$\vph_{u,A}(x)=\psi_{u,A}(\log|x_1|,\ldots,\log|x_n|)$, where
$$\psi_{u,A}(t)=\sup\,\{ \langle b,t\rangle:\: {b\in H_{u,A}}\},
\quad t\in \Rnm,$$ \beq H_{u,A}=\bigcap_{a\in A} \{b\in\Rnp:\:
\sum_k\frac{b_k}{a_k}\ge \nu_u(0,a)\}.\eeq Thus, if the only
information on $u$ is that it is locally bounded on
$D\setminus\{0\}$ and the values $\nu_u(0,a)$ for $a\in A$, then its
residual Monge--Amp\`ere mass at $0$ can be estimated as \beq
(dd^cu)^n(0)\ge (dd^ch^P_{u,0,D})^n(0)= (dd^c\vph_{u,A})^n(0) = n!\,
Vol(H_{u,A}).\eeq

(3) If $P$ consists of all the directional weights $\phi_{a,\zeta}$
(\ref{eq:phax}), then $h^P_{u,\zeta,D}=G_{\Psi,\zeta,D}$, the Green
function (\ref{eq:LNind}) with respect to the indicator
$\Psi=\Psi_{u,\zeta}$ (\ref{eq:lind}) of the function $u$ at
$\zeta$.

(4) Let $u=\log|z_1|$ in the unit polydisc $\D^n$ and let $P$ be the
collection of the weights $\phi_{a,0}$ (\ref{eq:phax}) with $a_1=1$.
Then the type of $u$ with respect to any $\phi_{a,0}\in P$ equals
$1$ and thus, $v\le \phi_{a,0}$ for every $v\in \M_{u,0,\D^n}$ and
any such direction $a$. Therefore $h^P_{u,0,\D^n}=u$. (For a more
general statement, see Theorem~\ref{theo:grf1}.)
\end{examples}

\begin{remarks} (1) As follows from Example 4, the functions
$h^P_{u,\zeta,D}$ need not be locally bounded outside $\zeta$.

(2) The same example shows that the condition of strong
plurisubharmonic barrier cannot be replaced by hyperconvexity in
general. However this can be done if $u\in
L_{loc}^\infty(\obl\setminus\{\zeta\})$, see
Proposition~\ref{prop:grf2}.
\end{remarks}

\subsection{Complete greenifications}

Another natural extremal function determined by the singularity of
$u$ can be defined as follows. Let $u\in PSH(\obl)$ be such that
$u(\zeta)=-\infty$ for some $\zeta\in\obl$. Given a domain
$D\subset\obl$, $\zeta\in D$, consider the class \beq
\F_{u,\zeta,D}=\{v\in PSH^-(D):\: v(z)\le u(z)+O(1),\
 z\to\zeta\},\eeq then the upper regularization of its upper envelope
is a plurisubharmonic function in $D$; we will call it the {\it
complete greenification} of $u$ at $\zeta$ and denote by
$g_{u,\zeta,D}$: \beq g_{u,\zeta,D}(x) = \limsup_{y\to
x}\sup\{v(y):\: v\in \F_{u,\zeta,D}\}.\eeq If $\vph\in MW_\zeta$,
then $g_{\vph,\zeta,D}=G_{\vph,D}$, the Green--Zahariuta function
(\ref{eq:gzf1}).

\medskip
It follows from the definition that if $u$ is bounded above on
$D$, then \beq\label{eq:uhu} u\le g_{u,\zeta,D}+\sup_Du.\eeq It is
easy to see that the function $g_{u,\zeta,D}$ need not belong to
the class $\F_{u,\zeta,D}$ (take, for example,
$u=-|\log|z||^{1/2}$, then $g_{u,\zeta,D}\equiv 0$).

\begin{proposition}\label{prop:grf0} If a plurisubharmonic function
$u$ is bounded above on $D\ni\zeta$, then
\begin{enumerate}
\item[(i)] $g_{u,\zeta,D}$ is maximal on any open $\omega\subset
D$ such that $g_{u,\zeta,D}\in L_{loc}^\infty(\omega)$;
\item[(ii)] $\nuph=\nu(g_{u,\zeta,D},\vph)$ for any continuous weight
$\vph$ with $\vph^{-1}(-\infty)=\zeta$; \item[(iii)]
$\suph=\sigma(g_{u,\zeta,D},\vph)$ for any $\vph\in MW_\zeta$;
\item[(iv)] if $D$ has a strong plurisubharmonic barrier at a point
$z\in
\partial D$ and if $u$ is bounded below near $z$, then $\lim_{x\to
z}g_{u,\zeta,D}(x)=0$.
\end{enumerate}
\end{proposition}

\begin{proof} Take a sequence of pseudoconvex domains $D_j$
such that $D_{j+1}\Subset D_j\Subset D$, $\cap_j D_j=\{\zeta\}$, and
let \beq\label{eq:grf1} u_j(x)=\sup\{v(x):\: v\in PSH^-(D), \ v\le u
-\sup_D \,u\quad {\rm in\ }D_j\},\quad x\in D.\eeq Since its upper
regularization $u_j^*$ belongs to $PSH^-(D)$ and coincides with $u
-\sup_D u$ in $D_j$, the function $u_j=u_j^*\in \F_{u,\zeta,D}$ and
 is maximal on $D\setminus D_j$. When
$j\to\infty$, the functions $u_j$ increase to a function $v$ such
that $v^*\in PSH^-(D)$ and is maximal where it is locally bounded.
Evidently, $g_{u,\zeta,D}\ge v^*$.

By the Choquet lemma, there exists a sequence
$w_k\in\F_{u,\zeta,D}$ that increases to $w$ such that
$w^*=g_{u,\zeta,D}$. Take any $\epsilon>0$, then for each $k$
there exists $j=j(k)$ such that $w_k\le (1-\epsilon)u_j$ on $D_j$.
Therefore $w_k\le (1-\epsilon)u_j\le (1-\epsilon)g_{u,\zeta,D}$ in
$D$, which gives $g_{u,\zeta,D}\le (1-\epsilon)v^*$ for all
$\epsilon>0$ and thus, $g_{u,\zeta,D}= v^*$, and (i) follows now
from the maximality of $v^*$.

To prove (ii), consider again the functions $u_j$ (\ref{eq:grf1}),
then $\nu(u_j,\vph)=\nuph$ for any $\vph$. By
Theorem~\ref{theo:Dsc}, \beq\nu(g_{u,\zeta,D},\vph)\ge
\limsup_{j\to\infty}\, \nu(u_j,\vph)=\nuph,\eeq while the reverse
inequality follows from (\ref{eq:uhu}) by Theorem~\ref{theo:CT2}.

Similar arguments (but now using Propositions~\ref{prop:converg} and
\ref{prop:rel} (v) instead of Theorems~\ref{theo:Dsc} and
\ref{theo:CT2}) prove (iii). Finally, (iv) can be proved exactly as
assertion (v) of Proposition~\ref{prop:genind0}.
\end{proof}

\medskip
More can be said if $u$ is locally bounded outside $\zeta$. Note
that then it can be extended (from a neighbourhood of $\zeta$) to
a plurisubharmonic function in the whole space, and none of its
greenifications at $\zeta$ depend on the choice of the extension.

\begin{proposition}\label{prop:grf2} If $D$ is a bounded
hyperconvex domain, $u\in PSH(D)\cap
L_{loc}^{\infty}(D\setminus\zeta)$, then $$\lim_{x\to z}
g_{u,\zeta,D}(x)=0,\quad z\in\partial D,$$
\beq\label{eq:gfr3}(dd^cu)^n({\zeta})=
(dd^cg_{u,\zeta,D})^n({\zeta}).\eeq
\end{proposition}

\begin{proof} The first statement follows exactly as in the case of
the pluricomplex Green function with logarithmic singularity.

To prove (\ref{eq:gfr3}), observe first that relation (\ref{eq:uhu})
implies $(dd^cu)^n({\zeta})\ge (dd^cg_{u,\zeta,D})^n({\zeta})$. On
the other hand, the functions $u_j$ (\ref{eq:grf1}) belong to the
Cegrell class $\F$ \cite{C2} and increase a.e. to $g_{u,\zeta,D}$.
By Theorem~5.4 of \cite{C2}, $(dd^cu_j)^n\to (dd^cg_{u,\zeta,D})^n$.
Therefore, $(dd^cg_{u,\zeta,D})^n({\zeta})\ge \limsup_{j\to\infty}
(dd^cu_j)^n({\zeta})=(dd^cu)^n({\zeta})$, which completes the proof.
\end{proof}


\begin{remark} In spite of the relations in
Proposition~\ref{prop:grf0} and (\ref{eq:gfr3}), some important
information on the singularity can be lost when passing to the
function $g_{u,\zeta,D}$. For example, if we take
$u(z)=\max\{\log|z_1|,-|\log|z_2||^{\frac 12}\}$, then $g_{u,0,D}=0$
for any $D\subset\C^2$ containing $0$, while $\nu(u,\log|z_2|)=1$
(note that the function $\log|z_2|$ is semiexhaustive on the support
of $dd^cu$).
\end{remark}

\section{Greenifications and Green functions}\label{sec:gr&gr}

Now we turn to relations between the extremal functions considered
above. We will write $\M_{u,\zeta,D}^S$ and $h_{u,\zeta,D}^S$ if
$P=SMW_\zeta$ (\ref{eq:SMW}), and $\M_{u,\zeta,D}^L$ and
$h_{u,\zeta,D}^L$ if $P=LMW_\zeta$ (\ref{eq:FMW}); in view of
(\ref{eq:sb}), $h^L_{u,\zeta,D}\ge h^S_{u,\zeta,D}$.

According to Proposition~\ref{prop:grf0}, we have
$\sigma(g_{u,\zeta,D},\vph)=\suph$ for all $\vph\in MW_\zeta$, so
\beq\label{eq:uhu1} g_{u,\zeta,D}\le h_{u,\zeta,D}\le
h^P_{u,\zeta,D}\eeq for any $u\in PSH(D)$ and $P\subset MW_\zeta$.

By Proposition~\ref{prop:ubind0}, the condition $g_{u,\zeta,D}\equiv
0$ implies $h_{u,\zeta,D}^P\equiv 0$ for every $P\subseteq
MW_\zeta$. So let us assume $g_{u,\zeta,D}\not\equiv 0$.

When $\vph\in MW_\zeta$, the function $g_{\vph,\zeta,D}$ is the
Green--Zahariuta function for the singularity $\vph$ in $D$; by
Proposition~\ref{prop:ubind0}, the same is true for
$h_{\vph,\zeta,D}$, so $g_{\vph,\zeta,D} = h_{\vph,\zeta,D}$. More
generally, we have the following simple

\begin{proposition}\label{prop:geqh}
Let $u\in PSH_\zeta$ be locally bounded outside $\zeta$, then
$g_{u,\zeta,D}= h_{u,\zeta,D}$. If, in addition, $\nu_u(\zeta)>0$,
then $g_{u,\zeta,D}= h_{u,\zeta,D}^S$.
\end{proposition}

\begin{proof} The equalities result from Proposition~\ref{prop:ubind0}
and (\ref{eq:uhu}) because $g_{u,\zeta,D}$ belongs to $MW_\zeta$
and, in case of $\nu_u(\zeta)>0$, to $SMW_\zeta$.
\end{proof}

\begin{remark} We do not know if $g_{u,\zeta,D}= h_{u,\zeta,D}^S$
when $\nu_u(\zeta)=0$. The condition $\nu_u(\zeta)=0$ implies, by
(\ref{eq:lelbound}), $h_{u,\zeta,D}^S\equiv 0$. As follows from
Theorem~\ref{prop:zer1} below, for functions $u$ locally bounded
outside $\zeta$ the relation $g_{u,\zeta,D}= h_{u,\zeta,D}^S$ is
thus equivalent to $(dd^cu)^n({\zeta})= 0$, and we are facing the
problem of existence of plurisubharmonic functions with zero Lelong
number and positive Monge-Amp\`ere mass. It can be reformulated as
follows: is it true that $g_{u,\zeta,D}\equiv 0$ if $u$ is locally
bounded outside $\zeta$ and $\nu_u(\zeta)=0$? Equivalently: is it
true that $SMW_\zeta=MW_\zeta\,$?
\end{remark}

To study the situation with the type-greenifications with respect to
arbitrary subsets $P$ of $MW_\zeta$, we need the following result on
"incommensurability" of Green functions.

\begin{lemma}\label{lem:zi} Let $D$ be a bounded hyperconvex
domain and let $v, w\in PSH(D)\cap
L_{loc}^{\infty}(D\setminus\zeta)$ be two solutions of the Dirichlet
problem $(dd^cu)^n=\tau\delta_\zeta$, $u|_{\partial D}=0$ with some
$\tau>0$. If $v\ge w$ in $D$, then $v\equiv w$.
\end{lemma}

\begin{proof} We use an idea from the proof of
Theorem~3.3 in \cite{Ze}. Choose $R>0$ such that
$\rho(x)=|x|^2-R^2<0$ in $D$. Given $\epsilon>0$, consider the
function $u_\epsilon=\max\{v+\epsilon\rho,w\}$. Since $u_\epsilon=w$
near $\partial D$, we have
\beq\label{eq:zer3}\int_D(dd^cu_\epsilon)^n=\int_D(dd^cw)^n=\tau.\eeq
On the other hand, $u_\epsilon\le v$ and thus, by
Theorem~\ref{theo:CT2}, \beq(dd^cu_\epsilon)^n({\zeta})\ge
(dd^cv)^n({\zeta})=\tau.\eeq Therefore (\ref{eq:zer3}) implies
$(dd^cu_\epsilon)^n=0$ on $D\setminus\{\zeta\}$. The functions
$v+\epsilon\rho$ and $w$ are locally bounded outside $\zeta$, so
\beq(dd^cu_\epsilon)^n\ge \chi_1(dd^c (v+\epsilon\rho))^n + \chi_2
(dd^cw)^n\eeq on $D\setminus\{\zeta\}$ (\cite{D4},
Proposition~11.9), where $\chi_1$ and $\chi_2$ are the
characteristic functions of the sets $E_1=\{w\le
v+\epsilon\rho\}\setminus\{\zeta\}$ and $E_2=\{w>
v+\epsilon\rho\}\setminus\{\zeta\}$, respectively. This gives
\beq\epsilon^n \int_{E_1}(dd^c\rho)^n\le \chi_1(dd^c
(v+\epsilon\rho))^n=0.\eeq Hence, the set $\{w < v\}$ has zero
Lebesgue measure, which proves the claim.
\end{proof}

\begin{theorem}\label{prop:zer1} Let $D$ be bounded and hyperconvex,
$u\in PSH(D)\cap L_{loc}^{\infty}(D\setminus\zeta)$, and let
$P\subseteq MW_\zeta$. Then $g_{u,\zeta,D}=h^P_{u,\zeta,D}$ if and
only if \beq\label{eq:zer2}(dd^cu)^n({\zeta})=
(dd^ch^P_{u,\zeta,D})^n({\zeta}).\eeq
\end{theorem}

\begin{proof} If $h^P_{u,\zeta,D}=g_{u,\zeta,D}$, then
(\ref{eq:zer2}) follows from Proposition~\ref{prop:grf2}. The
reverse implication follows from Lemma~\ref{lem:zi} (by
Proposition~\ref{prop:grf2} and (\ref{eq:uhu1}), the functions
$v=h^P_{u,\zeta,D}$ and $w=g_{u,\zeta,D}$ satisfy the conditions
of the lemma).
\end{proof}

\medskip
This can be applied to evaluation of the multiplicity of an
equidimensional holomorphic mapping by means of its Newton
polyhedron. Let $\zeta$ be an isolated zero of an equidimensional
holomorphic mapping $f$. Denote by $\Gamma_+(f,\zeta)$ the Newton
polyhedron of $f$ at $\zeta$, i.e., the convex hull of the set
$E_\zeta+\Rnp$, where $E_\zeta\subset\Znp$ is the collection of the
exponents in the Taylor expansions of the components of $f$ about
$\zeta$, and let $N_\zeta$ denote the Newton number of $f$ at
$\zeta$, i.e., $N_\zeta=n!\,{\rm Vol}(\Rnp\setminus
\Gamma_+(f,\zeta))$. Kouchnirenko's theorem \cite{Ku1} states that
the multiplicity $m_\zeta$ of $f$ at $\zeta$ can be estimated from
below by the Newton number, \beq\label{eq:NN} m_\zeta\ge
N_\zeta,\eeq with an equality under certain non-degeneracy
conditions. An application of Theorem~\ref{prop:zer1} gives a
necessary and sufficient condition for the equality to hold.

As was shown in \cite{R} and \cite{R4},
$N_\zeta=(dd^c\Psi_{u,\zeta})^n(0)$, where $\Psi_{u,\zeta}$ is the
indicator (\ref{eq:lind}) of the plurisubharmonic function
$u=\log|f|$ at $\zeta$. Let $P$ consist of all the directional
weights $\phi_{a,\zeta}$, $a\in\Rnp$, and let $D$ be a ball around
$\zeta$. Then (see Example~3 after Proposition~\ref{prop:ubind0})
the function $h^P_{u,\zeta,D}$ coincides with the Green function
(\ref{eq:LNind}) with respect to the indicator $\Psi_{u,\zeta}$,
which in turn equals the Green--Zahariuta function $G_{\vph,D}$
(\ref{eq:gzf}) for the singularity
$\vph(x)=\Psi_{u,\zeta}(x-\zeta)$. Therefore
$h^P_{u,\zeta,D}=g_{\vph,\zeta,D}$ and \beq N_\zeta=
(dd^ch^P_{u,\zeta,D})^n(\zeta).\eeq Since
$m_\zeta=(dd^cu)^n(\zeta)$, the equality $m_\zeta = N_\zeta$ is
equivalent to (\ref{eq:zer2}) and thus, by Theorem~\ref{prop:zer1},
to $g_{u,\zeta,D}=g_{\vph,\zeta,D}$. Finally, as $u,\vph\in
MW_\zeta$, we have $u=g_{u,\zeta,D}+O(1)$ and
$\vph=g_{\vph,\zeta,D}+O(1)$ near $\zeta$, which gives
$u=\vph+O(1)$. We have just proved the following

\begin{corollary}\label{cor:newt} The multiplicity of an
isolated zero $\zeta$ of an equidimensional holomorphic mapping $f$
equals its Newton number at $\zeta$ if and only if
$\log|f(x)|=\Psi(x-\zeta)+O(1)$ as $x\to\zeta$, where
$\Psi=\Psi_{\log|f|,\zeta}$ is the indicator (\ref{eq:lind}) of the
function $\log|f|$ at $\zeta$.
\end{corollary}

\medskip
The situation with non-isolated singularities looks more
complicated. Observe, for example, that $h_{u,\zeta,D}^P$ is maximal
on the whole $D\setminus\{\zeta\}$, however we do not know if the
same is true for $g_{u,\zeta,D}$.

We can handle the situation in the case of analytic singularities.
In this section, we prove the equality
$g_{u,\zeta,D}=h_{u,\zeta,D}^L$ for $u=\log|f|$, where $f:D\to\C$ is
a holomorphic function; we recall that $h_{u,\zeta,D}^L$ is the
type-greenification with respect to the class $LMW_\zeta$
(\ref{eq:FMW}). In this case, the greenifications coincide with the
Green function $G_{A,D}$ (\ref{eq:LSG1}) in the sense of
L\'arusson--Sigurdsson. Note that the function $G_{A,D}$ is defined
as the upper envelope of functions $u$ with $\nu_u(a)\ge \nu_A(a)$
for {\sl all} $a\in |A|$. It turns out that one can consider only
one point (or finitely many ones) from $|A|$, but then an infinite
set of weights should be used. The case of mappings $f:D\to\C^p$,
$p>1$, will be considered in Section~\ref{sec:Hoelder}.

\begin{theorem}\label{theo:grf1} Let $u=\log|f|$, where $f$
is a holomorphic function on $\Omega$. Given $\zeta\in f^{-1}(0)$,
let $f=ab$ with $a=a_1^{m_1}\ldots a_k^{m_k}$ such that $a_j$ are
irreducible factors of $f$ vanishing at $\zeta$ and $b(\zeta)\neq
0$. Then for any hyperconvex domain $D\Subset\obl$ that contains
$\zeta$,
\beq\label{eq:sldiv}h_{u,\zeta,D}^L=g_{u,\zeta,D}=G_{A_\zeta,D},\eeq
where $G_{A_\zeta,D}$ is the Green function (\ref{eq:LSG1}) for the
divisor $A_\zeta$ of the function $a$. Moreover, there exists a
sequence $P$ of continuous weights $\vph_j\in LMW_\zeta$ such that
the Green--Zahariuta functions $G_{\vph_j,D}$ decrease to
$h^P_{u,\zeta,D}=G_{A_\zeta,D}$.
\end{theorem}

\begin{proof}
Let $\F_{A_\zeta,D}$ be the class defined by (\ref{eq:LSG2}) for
$A=A_\zeta$, then $\F_{A_\zeta,D}\subset\F_{u,\zeta,D}\subset
\M_{u,\zeta,D}^L$, so \beq G_{A_\zeta,D}\le g_{u,\zeta,D}\le
h_{u,\zeta,D}^L.\eeq

Choose a sequence of domains $D_j$ such that $D_{j+1}\Subset
D_j\Subset D$, $\cap_j D_j=\{\zeta\}$, and $fa^{-1}$ does not vanish
in $D_1$, then the functions $u_j$ defined by (\ref{eq:grf1})
satisfy \beq\label{eq:lsgr}u_j\le \log|a|+C_j\quad {\rm in}\
D_j.\eeq
 By Siu's theorem, the set
$\{x\in D: \nu(u_j,x)\ge m_s\}$ is analytic; by (\ref{eq:lsgr}), it
contains the support $|A_s|$ of the divisor $A_s$ of $a_s$.
Therefore, $\nu(u_j,x)\ge\nu(\log|a|,x)$ for all $x\in |A_\zeta|$.
Since $u_j$ converge to $g_{u,\zeta,D}$, this implies
$\nu(g_{u,\zeta,D},x)\ge\nu(\log|a|,x)$ and thus,
$g_{u,\zeta,D}\in\F_{A_\zeta,D}$. This proves the second equality in
(\ref{eq:sldiv}).

Let $f_2,\ldots,f_n$ be holomorphic functions in a neighbourhood
$\omega\Subset D$ of $\zeta$ such that $\zeta$ is the only point of
the zero set of the mapping $(f_1,\ldots,f_n)$ in $\omega$ and
$\Omega_0=\{z\in\omega:\: |f_k|<1,\ 1\le k\le n\}\Subset \omega$,
where $f_1=2a(\sup_\omega|a|)^{-1}$. Denote
\beq\vph_j:=\sup\{\log|f_1|,j\log|f_k|:\:2\le k\le n\},\quad
j\in\Z_+.\eeq We have $\vph_j\in LMW_\zeta(\Omega_0)$ and $\vph_j=0$
on $\partial\Omega_0$, so $\vph_j=G_{\vph_j,\Omega_0}$, the
Green--Zahariuta function for the singularity $\vph_j$ in
$\Omega_0$. Since $h_{u,\zeta,D}^L<0$ in $\Omega_0$ and
$\sigma(h_{u,\zeta,D}^L,\vph_j)=\sigma(u,\vph_j)=1$, we have
$h_{u,\zeta,D}^L\le \vph_j$ in $\Omega_0$ for each $j$. Therefore
$h_{u,\zeta,D}^L\le \log|f_1|=\lim_{j\to\infty}\vph_j$ in
$\Omega_0$. This implies $h_{u,\zeta,D}^L\in \F_{u,\zeta,D}$ and
thus, $h_{u,\zeta,D}^L=g_{u,\zeta,D}$.

Finally, the Green--Zahariuta functions $G_{\vph_j,D}$ dominate
$h_{u,\zeta,D}^L$ and decrease to some function $v\in PSH^-(D)$.
Since
$\sigma(v,\vph_k)\ge\limsup_{j\to\infty}\sigma(\vph_j,\vph_k)=1$, we
get $v\le h_{u,\zeta,D}^L$, which completes the proof.
\end{proof}

\bigskip
One can also consider the greenifications with respect to arbitrary
finite sets ${\cal Z}\subset D$, \beq\label{eq:hsev} h_{u,{\cal
Z},D}^P(x)=\sup\,\{v(x):\:v\in \M_{u,\zeta,D}^{P(\zeta)},\ \zeta\in
{\cal Z}\}\eeq and \beq\label{eq:gsev} g_{u,{\cal Z},D}(x) =
\limsup_{y\to x}\sup\{v(y):\: v\in \F_{u,\zeta,D},\ \zeta\in {\cal
Z}\}.\eeq They have properties similar to those of the functions
$h_{u,\zeta,D}^P$ and $g_{u,\zeta,D}$, stated in
Propositions~\ref{prop:genind0}--\ref{prop:grf2} (with obvious
modifications). In particular, if $\vph\in PSH(D)$ is such that
$\vph^{-1}(-\infty)={\cal Z}$ and  the restriction of $\vph$ to a
neighbourhood of $\zeta$ belongs to $MW_\zeta$ for each $\zeta\in
{\cal Z}$, then $ h_{\vph,{\cal Z},D}=g_{\vph,{\cal
Z},D}=G_{\vph,D}$, the Green--Zahariuta function with the
singularities defined by $\vph$. They are also related to the Green
functions (\ref{eq:LSG1}) as follows (cf. Theorem~\ref{theo:grf1}).

\begin{theorem}\label{theo:grf2} Let $A$ be the divisor of
a holomorphic function $f$ in $\Omega$ and let $u=\log|f|$. If a
finite subset ${\cal Z}$ of a hyperconvex domain $D\Subset\Omega$ is
such that each irreducible component of $|A|\cap D$ contains at
least one point of ${\cal Z}$, then $ h_{u,{\cal Z},D}^L=g_{u,{\cal
Z},D}=G_{A,D}$.
\end{theorem}


\section{Analyticity theorem}\label{sec:Hoelder}

We let $\obl$ be a  pseudoconvex domain in $\Cn$, and let $R:\obl\to
(-\infty,\infty]$ be a lower semicontinuous function on $\obl$. We
consider a continuous plurisubharmonic function
$\vph:\obl\times\obl\to [-\infty,\infty)$ such that:
\begin{enumerate} \item[(i)] $\vph(x,\zeta)<R(\zeta)$ on $\obl\times\obl$;
\item[(ii)] $\{x:\vph(x,\zeta)= -\infty\}=\{\zeta\}$; \item[(iii)] for any
$\zeta\in\obl$ and $r<R(\zeta)$ there exists a neighbourhood $U$ of
$\zeta$ such that the set $\{(x,y):\: \vph(x,y)<r,\: y\in U\}\Subset
\Omega\times \obl$; \item[(iv)] $(dd^c\vph)^n=0$ on $\{
\vph(x,\zeta)>-\infty\}$; \item[(v)] $e^{\vph(x,\zeta)}$ is H\"older
continuous in $\zeta$: \beq\label{eq:holder}\exists\beta>0:\
|e^{\vph(x,\zeta)}-e^{\vph(x,y)}|\le |\zeta-y|^\beta,\quad x,y,\zeta
\in\obl.\eeq
\end{enumerate}

It then follows that $\vph_\zeta(x):=\vph(x,\zeta)\in SMW_\zeta$
(\ref{eq:SMW}). Similarly to (\ref{eq:luvr}) and (\ref{eq:type}) we
introduce the function $\Lambda(u,\vphz,r)$ and the relative type
$\sigma(u,\vphz)$. By Theorem~6.8 of \cite{D3}, $\Lambda(u,\vphz,r)$
is plurisubharmonic on each connected component of the set $\{\zeta:
R(\zeta)>r\}$.

\medskip
As was shown (in a more general setting) by Demailly \cite{D1}, the
sets $\{\zeta:\:\nu(u,\vph_\zeta)\ge c\}$ are analytic for all
$c>0$. By an adaptation of Kiselman's and Demailly's proofs of Siu's
theorem, we prove its analogue for the relative types. Denote \beq
S_c(u,\vph,\Omega)=\{\zeta\in\obl:\: \sigma(u,\vphz)\ge c\}, \quad
c>0. \eeq

\begin{theorem}\label{theo:typean} Let $\obl$ be a
pseudoconvex domain in $\Cn$, a function $\vph(x,\zeta)$ satisfy the
above conditions (i)--(v), and $u\in PSH(\obl)$. Then
$S_c(u,\vph,\Omega)$ is an analytic subset of $\obl$ for each $c>0$.
\end{theorem}

\begin{proof}
By Theorem~6.11 of \cite{D3}, the function
$U(\zeta,\xi)=\Lambda(u,\vphz,{\rm Re}\,\xi)$ is plurisubharmonic in
$\{(\zeta,\xi)\in \obl\times\C:\: {\rm Re}\,\xi<R(\zeta)\}$. Fix a
pseudoconvex domain $D\Subset\obl$ and denote
$R_0=\inf\,\{R(\zeta):\: \zeta\in D\}>-\infty$. Given $a>0$, the
function $U(\zeta,\xi)-a\,{\rm Re}\,\xi$ is then plurisubharmonic
and independent of ${\rm Im}\,\xi$, so by Kiselman's minimum
principle \cite{Kis1}, the function \beq
U_a(\zeta)=\inf\{\Lambda(u,\vphz,r)-ar:\:r<R_0\} \eeq is
plurisubharmonic in $D$.

Let $\zeta\in D$. If $a>\sigma(u,\vphz)$, then $\Lambda(u,\vphz,r)>
ar$ for all $r\le r_0<\min\,\{R_0,0\}$. If $r_0<r<R_0$, then
$\Lambda(u,\vphz,r)- ar> \Lambda(u,\vph_y,r_0)- aR_0$. Therefore
$U_a(\zeta)>-\infty$.

Now let $a<\sigma(u,\vphz)$. H\"older continuity (\ref{eq:holder})
implies
\beq\Lambda(u,\vph_y,r)\le
\Lambda(u,\vphz,\log(e^r+|y-\zeta|^\beta))\le
\sigma(u,\vphz)\log(e^r+|y-\zeta|^\beta)+C\eeq in a neighbourhood of
$\zeta$. Denote $r_y=\beta\log|y-\zeta|$, then \beq\label{eq:holb1}
U_a(y)\le \Lambda(u,\vph_y,r_y)- ar_y \le
(\sigma(u,\vphz)-a)\beta\log|y-\zeta|+C_1\eeq near $\zeta$.

Finally, let $Z_{a,b}$ ($a,b>0$) be the set of points $\zeta\in D$
such that $\exp(-b^{-1}U_a)$ is not integrable in a neighbourhood of
$\zeta$. As follows from the H\"ormander--Bombieri--Skoda theorem,
the sets $Z_{a,b}$ are analytic.

If $\zeta\not\in S_c(u,\vph,\Omega)$, then the function $b^{-1}U_a$
with $\sigma(u,\vphz)<a<c$ is finite at $\zeta$ and so, by Skoda's
theorem, $\zeta\not\in Z_{a,b}$ for all $b>0$. If $\zeta\in
S_c(u,\vph,\Omega)$, $a<c$, and $b< (c-a)\beta (2n)^{-1}$, then
(\ref{eq:holb1}) implies $\zeta\in Z_{a,b}$. Thus,
$S_c(u,\vph,\Omega)$ coincides with the intersection of all the sets
$ Z_{a,b}$ with $a<c$ and $b<(c-a)\beta (2n)^{-1}$, and is therefore
analytic.
\end{proof}

\begin{remark} The result can be reformulated in the following way:
under the conditions of Theorem~\ref{theo:typean}, the set
$S(u,\vph,\Omega)=\{\zeta\in\Omega:\: u(x)\le\vph(x,\zeta) +O(1)\
{\rm as\ }x\to\zeta\}$ is analytic. As was noticed by the referee,
condition (iv) is actually necessary. Take, for example, the
function
$\vph(x,\zeta)=\max\{\log|x_1-\zeta_1|+\log|(x_1-\zeta_1)x_2|,
\log|x_2-\zeta_2|\}$ in $\C^2\times\C^2$; it has all the properties
except for (iv), while the set
$S(\log|x_1|,\vph,\C^2)=\{(0,\zeta_2):\zeta_2\neq 0\}$ is not
analytic.
\end{remark}

As an application, we present the following result on propagation of
plurisubharmonic singularities. We will say that a closed complex
space $A$ is a {\it locally complete intersection} if $|A|$ is of
pure codimension $p$ and the associated ideal sheaf $\I_A$ is
locally generated by $p$ holomorphic functions.

\begin{corollary}\label{cor:propag} Let a closed complex
space $A$ be locally complete intersection on a domain
$D\subset\Cn$, and let $\omega$ be an open subset of $D$ that
intersects each irreducible component of $|A|$. If a function $u\in
PSH(D)$ satisfies $u\le\log|f|+O(1)$ locally in $\omega$, where
$f_1,\ldots,f_p$ are local generators of $\I_A$, then it satisfies
this relation near every point of $|A|$.
\end{corollary}

\begin{proof} Denote by $Z_l$, $l=1,2,\ldots$, the irreducible
components of $|A|$, and \beq Z_l^*=(Reg\,Z_l)\setminus
\bigcup_{k\neq l}Z_k.\eeq

Every point $z\in Z_l^*$ has a neighbourhood $V$ and coordinates
$x=(x',x'')\in \C^p\times \C^{n-p}$, centered at $z$, such that
$V\cap |A|=V\cap Z_l^*= V\cap \{x'=0\}$ and $f_1,\ldots,f_p$ are
global generators of $\I_A$ on $V$. Let $U\Subset D$ be a
pseudoconvex neighbourhood of $z$ such that $U-U\subset V$. For any
$\zeta\in U$ and $N>0$, the function
\beq\vph^N(x,\zeta)=\max\,\{\log|f(x-(\zeta',0))|,N\log|x''-\zeta''|\}\eeq
satisfies the conditions of Theorem~\ref{theo:typean} with
$\Omega=U$. Therefore, $S_1(u,\vph^N,U)$ is an analytic subset of
$U$. Note that $S_1(\log|f|,\vph^N,U)=U\cap |A|$.

Let $\{U^j\}$ be a denumerable covering of $Z_l^*$ by such
neighbourhoods, and let $\{\vph^{j,N}\}$ be the corresponding
weights. We may assume $u\le\log|f|+O(1)$ near a point $\zeta_0\in
U^1\cap\omega\cap |A|$, so $\sigma(u,\vph^{1,N}_\zeta)\ge
\sigma(\log|f|,\vph^{1,N}_\zeta)= 1$ for every $\zeta\in
U^1\cap\omega\cap |A|$ and thus for every $\zeta\in U^1\cap |A|$.

Take any $\zeta\in U^1\cap Z_l^*$ and choose constants $a,b>0$ such
that \beq V_\zeta=\{x\in U^1:\: a|f(x)|<1, \
b|x''-\zeta''|<1\}\Subset U^1.\eeq Assuming $u\le C$ in $U^1$, the
relation $\sigma(u,\vph^{j,N}_\zeta)\ge 1$ implies $u\le g^{j,N}+C$
in $V_\zeta$, where \beq g^{j,N}(x)= \max\,\{\log|af(x)|,N\log
b|x''-\zeta''|\}\eeq is the Green--Zahariuta function in $V_\zeta$
for the singularity $\vph^{j,N}_\zeta$. Observe  that $g^{j,N}$
decrease to $\log|af|$ as $N\to\infty$, so $u\le\log|f|+C_1$ near
$\zeta$. Therefore we have extended the hypothesis of the theorem
from $\omega$ to $\omega\cup U^1$. By repeating the argument, we get
$u\le\log|f|+O(1)$ near every point of $Z_l^*$ (since it is
connected) and so, of $Reg\,|A |$. Finally, the bound near irregular
points of $|A|$ can be deduced by using Thie's theorem and the
equation $(dd^c\log|f|)^p=0$ outside $|A|$ (see, for example, the
proof of Lemma~4.2 in \cite{RSig2}).
\end{proof}

\begin{remarks} (1) In particular, if $F$ is a holomorphic function in
$D$ whose restriction to $\omega$ belongs to the integral closure
${\bar\I}_{A,\omega}$ of the ideal sheaf $\I_{A,\omega}$ on
$\omega$, then $F\in {\bar\I}_{A,D}$. This follows from
Corollary~\ref{cor:propag} because $F\in {\bar\I}_A$ if and only if
$\log|F|\le\log|f|+O(1)$.

(2) Corollary~\ref{cor:propag} fails when the complete intersection
assumption is removed. Take, for example, $f=(z_1^2,z_1z_2)$ in
$\C^2$, then $\log|z_1|\le\log|f(z)|+O(1)$ near every point
$(0,\xi)$ with $\xi\neq 0$, but not near the origin.
\end{remarks}

A consequence of Corollary~\ref{cor:propag} is the following
analogue to Theorems \ref{theo:grf1} and \ref{theo:grf2} for the
Green functions (\ref{eq:RS}) with singularities along complex
spaces.

\begin{corollary}\label{cor:grf3} Let a closed complex
space $B$ be locally complete intersection on a pseudoconvex domain
$\Omega\subset\Cn$, and let $A$ be the restriction of $B$ to a
hyperconvex domain $D\Subset \Omega$. Further, let $F_1,\ldots,F_m$
be global sections of $\I_B$ generating $\I_A$, and $u=\log|F|$. If
a finite subset ${\cal Z}$ of $D$ is such that each irreducible
component of $|A|$ contains at least one point of ${\cal Z}$, then
\beq\label{eq:sldiv3} h_{u,{\cal Z},D}^L=g_{u,{\cal
Z},D}=G_{A,D},\eeq where the functions $h_{u,{\cal Z},D}^L$ and
$g_{u,{\cal Z},D}$ are defined by (\ref{eq:hsev}) and
(\ref{eq:gsev}), $L$ denotes the collection of maximal weights with
finite \L ojasiewicz exponents (\ref{eq:uln}) at the points of $\cal
Z$, and $G_{A,D}$ is the Green function (\ref{eq:RS}) for the space
$A$.
\end{corollary}

\begin{proof} We use the same arguments as in the proof of
Theorems \ref{theo:grf1} and \ref{theo:grf2}, the only difference
being referring to Corollary~\ref{cor:propag} instead of Siu's
theorem.

By (\ref{eq:RS1}), we have $G_{A,D}\le g_{u,{\cal Z},D}$. Let ${\cal
Z}=\{\zeta_1,\ldots,\zeta_s\}$. For $k=1,\ldots,s$ take a sequence
of pseudoconvex domains $D_{k,j}$ such that $D_{k,j+1}\Subset
D_{k,j}\Subset D$, $\cap_j D_{k,j}=\{\zeta_k\}$, and denote
$D_j=\cup_k D_{k,j}$. As in the proof of
Proposition~\ref{prop:grf0}, the functions $u_j$ defined by
(\ref{eq:grf1}) with such a choice of $D_j$ are plurisubharmonic in
$D$, the sequence increases to $g_{u,{\cal Z},D}$ a.e. in $D$, and
$u_j\le \log|F|+C_j$ near each $\zeta_k$. Since
$\log|F|=\log|f|+O(1)$, where $f_1,\ldots,f_p$ are local generators
of $\I_A$, $p=\codim |A|$, this gives, by
Corollary~\ref{cor:propag}, the relations $u_j\le\log|f|+O(1)$
locally near $|A|$ and thus, $g_{u,{\cal Z},D}\le G_{A,D}$. This
proves the second equality in (\ref{eq:sldiv3}).

The rest is proved as in Theorem~\ref{theo:grf1}.
\end{proof}

\bigskip
{\small {\it Acknowledgement.} The author is grateful to the referee
for valuable remarks and suggestions that have helped much in
improving the presentation.}

\vskip1cm

Tek/Nat, University of Stavanger, 4036 Stavanger, Norway

\vskip0.1cm

{\sc E-mail}: alexander.rashkovskii@uis.no

\end{document}